\documentclass[11pt, a4paper]{article}
\usepackage{amsmath} \usepackage{euscript}
\usepackage{mymatrix2}\let\mymatrixII\mymatrix
\usepackage{mymatrix}
\usepackage{longtable}
\usepackage{amssymb}
\begin{document}

\newcounter{bnomer} \newcounter{snomer}
\newcounter{bsnomer}
\setcounter{bnomer}{0}
\renewcommand{\thesnomer}{\thebnomer.\arabic{snomer}}
\renewcommand{\thebsnomer}{\thebnomer.\arabic{bsnomer}}
\renewcommand{\refname}{\begin{center}\large{\textbf{References}}\end{center}}

\newcommand{\sect}[1]{%
\setcounter{snomer}{0}\setcounter{bsnomer}{0}
\refstepcounter{bnomer}
\par\bigskip\begin{center}\large{\textbf{\arabic{bnomer}. {#1}}}\end{center}}
\newcommand{\sst}{%
\refstepcounter{bsnomer}
\par\bigskip\textbf{\arabic{bnomer}.\arabic{bsnomer}. }}
\newcommand{\defi}[1]{%
\refstepcounter{snomer}
\par\textbf{Definition \arabic{bnomer}.\arabic{snomer}. }{#1}\par}
\newcommand{\theo}[2]{%
\refstepcounter{snomer}
\par\textbf{Theorem \arabic{bnomer}.\arabic{snomer}. }{#2} {\emph{#1}}\hspace{\fill}$\square$\par}
\newcommand{\mtheo}[1]{%
\refstepcounter{snomer}
\par\textbf{Theorem \arabic{bnomer}.\arabic{snomer}. }{\emph{#1}}\par}
\newcommand{\theobp}[2]{%
\refstepcounter{snomer}
\par\textbf{Theorem \arabic{bnomer}.\arabic{snomer}. }{#2} {\emph{#1}}\par}
\newcommand{\theop}[2]{%
\refstepcounter{snomer}
\par\textbf{Theorem \arabic{bnomer}.\arabic{snomer}. }{\emph{#1}}
\par\textsc{Proof}. {#2}\hspace{\fill}$\square$\par}
\newcommand{\theosp}[2]{%
\refstepcounter{snomer}
\par\textbf{Theorem \arabic{bnomer}.\arabic{snomer}. }{\emph{#1}}
\par\textbf{Sketch of the proof}. {#2}\hspace{\fill}$\square$\par}
\newcommand{\exam}[1]{%
\refstepcounter{snomer}
\par\textbf{Example \arabic{bnomer}.\arabic{snomer}. }{#1}\par}
\newcommand{\deno}[1]{%
\refstepcounter{snomer}
\par\textbf{Definition \arabic{bnomer}.\arabic{snomer}. }{#1}\par}
\newcommand{\post}[1]{%
\refstepcounter{snomer}
\par\textbf{Proposition \arabic{bnomer}.\arabic{snomer}. }{\emph{#1}}\hspace{\fill}$\square$\par}
\newcommand{\postp}[2]{%
\refstepcounter{snomer}
\par\textbf{Proposition \arabic{bnomer}.\arabic{snomer}. }{\emph{#1}}
\par\textsc{Proof}. {#2}\hspace{\fill}$\square$\par}
\newcommand{\lemm}[1]{%
\refstepcounter{snomer}
\par\textbf{Lemma \arabic{bnomer}.\arabic{snomer}. }{\emph{#1}}\hspace{\fill}$\square$\par}
\newcommand{\lemmp}[2]{%
\refstepcounter{snomer}
\par\textbf{Lemma \arabic{bnomer}.\arabic{snomer}. }{\emph{#1}}
\par\textsc{Proof}. {#2}\hspace{\fill}$\square$\par}
\newcommand{\coro}[1]{%
\refstepcounter{snomer}
\par\textbf{Corollary \arabic{bnomer}.\arabic{snomer}. }{\emph{#1}}\hspace{\fill}$\square$\par}
\newcommand{\mcoro}[1]{%
\refstepcounter{snomer}
\par\textbf{Corollary \arabic{bnomer}.\arabic{snomer}. }{\emph{#1}}\par}
\newcommand{\corop}[2]{%
\refstepcounter{snomer}
\par\textbf{Corollary \arabic{bnomer}.\arabic{snomer}. }{\emph{#1}}
\par\textsc{Proof}. {#2}\hspace{\fill}$\square$\par}
\newcommand{\nota}[1]{%
\refstepcounter{snomer}
\par\textbf{Remark \arabic{bnomer}.\arabic{snomer}. }{#1}\par}
\newcommand{\propp}[2]{%
\refstepcounter{snomer}
\par\textbf{Proposition \arabic{bnomer}.\arabic{snomer}. }{\emph{#1}}
\par\textsc{Proof}. {#2}\hspace{\fill}$\square$\par}
\newcommand{\hypo}[1]{%
\refstepcounter{snomer}
\par\textbf{Conjecture \arabic{bnomer}.\arabic{snomer}. }{\emph{#1}}\par}

\newcommand{\Ind}[3]{%
\mathrm{Ind}_{#1}^{#2}{#3}}
\newcommand{\Res}[3]{%
\mathrm{Res}_{#1}^{#2}{#3}}
\newcommand{\epsi}{\epsilon}
\newcommand{\tri}{\triangleleft}
\newcommand{\Supp}[1]{%
\mathrm{Supp}(#1)}

\newcommand{\reg}{\mathrm{reg}}
\newcommand{\sreg}{\mathrm{sreg}}
\newcommand{\codim}{\mathrm{codim}\,}
\newcommand{\chara}{\mathrm{char}\,}
\newcommand{\rk}{\mathrm{rk}\,}
\newcommand{\id}{\mathrm{id}}
\newcommand{\Ad}{\mathrm{Ad}}
\newcommand{\col}{\mathrm{col}}
\newcommand{\row}{\mathrm{row}}
\newcommand{\low}{\mathrm{low}}
\newcommand{\pho}{\hphantom{\quad}\vphantom{\mid}}
\newcommand{\wt}{\widetilde}
\newcommand{\wh}{\widehat}
\newcommand{\ad}[1]{\mathrm{ad}_{#1}}
\newcommand{\tr}{\mathrm{tr}\,}
\newcommand{\GL}{\mathrm{GL}}
\newcommand{\Sp}{\mathrm{Sp}}
\newcommand{\Mat}{\mathrm{Mat}}
\newcommand{\Stab}{\mathrm{Stab}}

\newcommand{\vfi}{\varphi}
\newcommand{\teta}{\vartheta}
\newcommand{\lee}{\leqslant}
\newcommand{\gee}{\geqslant}
\newcommand{\Fp}{\mathbb{F}}
\newcommand{\Rp}{\mathbb{R}}
\newcommand{\Zp}{\mathbb{Z}}
\newcommand{\Cp}{\mathbb{C}}
\newcommand{\ut}{\mathfrak{u}}
\newcommand{\at}{\mathfrak{a}}
\newcommand{\nt}{\mathfrak{n}}
\newcommand{\spt}{\mathfrak{sp}}
\newcommand{\rt}{\mathfrak{r}}
\newcommand{\rad}{\mathfrak{rad}}
\newcommand{\bt}{\mathfrak{b}}
\newcommand{\gt}{\mathfrak{g}}
\newcommand{\vt}{\mathfrak{v}}
\newcommand{\pt}{\mathfrak{p}}
\newcommand{\Po}{\EuScript{P}}
\newcommand{\Uo}{\EuScript{U}}
\newcommand{\Fo}{\EuScript{F}}
\newcommand{\Do}{\EuScript{D}}
\newcommand{\Eo}{\EuScript{E}}
\newcommand{\Iu}{\mathcal{I}}
\newcommand{\Mo}{\mathcal{M}}
\newcommand{\Nu}{\mathcal{N}}
\newcommand{\Ro}{\mathcal{R}}
\newcommand{\Co}{\mathcal{C}}
\newcommand{\Lo}{\mathcal{L}}
\newcommand{\Ou}{\mathcal{O}}
\newcommand{\Au}{\mathcal{A}}
\newcommand{\Vu}{\mathcal{V}}
\newcommand{\Bu}{\mathcal{B}}
\newcommand{\Sy}{\mathcal{Z}}
\newcommand{\Sb}{\mathcal{F}}
\newcommand{\Gr}{\mathcal{G}}

\author{Mikhail V. Ignatyev\thanks{The research was supported by RFBR grant no. 11--01--90703--mob${}_-$st.}}
\date{\small Department of Algebra and Geometry\\
Samara State University\\
Samara, 443011, Ak.\ Pavlova, 1, Russia\\
\texttt{mihail.ignatev@gmail.com}}
\title{\Large{The Bruhat--Chevalley order on involutions of the hyperoctahedral group and~combinatorics of $B$-orbit closures}} \maketitle

\vspace{-1cm}\sect{Introduction and the main result}

\sst Let $G$ be a complex reductive algebraic group, $B$ a Borel subgroup of $G$, $\Phi$ the root system of~$G$ and $W=W(\Phi)$ the Weyl group of $\Phi$. It is well-known that the Bruhat--Chevalley order on~$W$ encodes the cell decomposition of the flag variety $G/B$ (see, e.g., \cite{BileyLakshmibai1}). Denote by $\Iu(W)$ the poset of involutions in~$W$ (i.e., elements of~$W$ of order 2). In \cite{RichardsonSpringer1}, R. Richardson and T.~Springer showed that $\Iu(A_{2n})$ encodes the incidences among the closed $B$-orbits on the symmetric variety $\mathrm{SL}_{2n+1}(\Cp)/\mathrm{SO}_{2n+1}(\Cp)$. In~\cite{BagnoCherniavsky1}, E.~Bagno and Y.~{Cher\-nav\-sky} presented a~{geo\-met\-ri\-cal} {inter\-pre\-ta\-tion} of the poset $\Iu(A_n)$, considering the action of the Borel subgroup of $\GL_n(\Cp)$ on symmetric matrices by congruence. F.~Incitti studied the poset $\Iu(\Phi)$ from a purely combinatorial point of view for the case of classical root system~$\Phi$ (see \cite{Incitti1}, \cite{Incitti2}). In particular, he proved that this poset is graded, calculated the rank function and described the covering relation.

In \cite{Ignatev3}, we presented another geometrical interpretation of $\Iu(A_n)$ in terms of \emph{coadjoint} $B$-orbits. Precisely, let $G=\GL_n(\Cp)$, the general linear group, then $\Phi=A_{n-1}$ and $W=S_n$, the symmetric group on $n$ letters. Let $B=B_n$ be the group of invertible upper-triangular matrices. Denote by $U=U_n$ the unitriangular group, i.e., the group of upper-triangular matrices with $1$'s on the diagonal. (In fact, $U$ is the unipotent radical of $B$.) Let $\nt=\nt_n$ be the space of upper-triangular matrices with zeroes on the diagonal, and $\nt^*$ the dual space. Since $B$ acts on $\nt$ by conjugation, one can consider the dual action of $B$ on $\nt^*$. To each involution $\sigma\in\Iu(A_{n-1})$ one can assign the $B$-orbit $\Omega_{\sigma}\subseteq\nt^*$ (see\linebreak Subsection~\ref{sst:A_n} or \cite[Subsection 1.2]{Ignatev3} for precise definitions). By~\cite[Theorem 1.1]{Ignatev3} (see also Theorem~\ref{mtheo_A_n}), $\Omega_{\sigma}$ is contained in the Zariski closure of $\Omega_{\tau}$ if and only if $\sigma$ is less or equal to $\tau$ with respect to the Bruhat--Chevalley order. In some sense, these results are ``dual'' to A. Melnikov's results \cite{Melnikov1}, \cite{Melnikov2}, \cite{Melnikov3}.

In this paper, we find similar results for the case $\Phi=C_n$. Namely, let $G=\Sp_{2n}(\Cp)$, the symplectic group, then $\Phi=C_n$ and $W$ is the hyperoctahedral group. Let $U$ be the unipotent radical of $B$, $\nt$ its Lie algebra and $\nt^*$ the dual space. Since $\nt$ is invariant under the adjoint action of the Borel subgroup on its Lie algebra, one can consider the dual action of $B$ on $\nt^*$. To each involution $\sigma\in\Iu(C_n)$ one can assign the $B$-orbit $\Omega_{\sigma}\subseteq\nt^*$ (see Definition \ref{defi:B_orbit_assoc_invol}). The main result of the paper is as follows.

\medskip\mtheo{Let $\sigma$, $\tau$ be involutions in the Weyl group of $C_n$. The orbit $\Omega_{\sigma}$ is contained in the Zariski closure of $\Omega_{\tau}$ if and only if $\sigma$ is less or equal to $\tau$ with respect to the Bruhat--Chevalley order.\label{mtheo}}

\medskip The paper is organized as follows. In the remainder of this section, we briefly recall our results about $\Iu(A_{n-1})$ (see Subsection \ref{sst:A_n}) and give precise definitions for the case of $C_n$ (see Subsection~\ref{sst:C_n}). Section \ref{sect:proof_mtheo} is devoted to the proof of Theorem~\ref{mtheo}. Precisely, in Subsection \ref{sst:if_Bruhat_then_contained}, using Incitti's results, we show that if $\sigma,\tau\in\Iu(C_n)$, then $\sigma\leq_B\tau$ implies $\Omega_\sigma\subseteq\overline{\Omega}_{\tau}$, see Proposition \ref{prop:if_Bruhat_then_contained}. (Here $\leq_B$ denotes the Bruhat--Chevalley order and $\overline{Z}$ denotes the Zariski closure of a subset $Z\subseteq\nt^*$.) Then, in Subsection~\ref{sst:star_equiv_Bruhat}, we define a partial order $\leq^*$ on $\Iu(C_n)$ in combinatorial terms and, using\linebreak \cite[Theorem~1.10]{Ignatev3}, prove Proposition \ref{prop:star_equiv_Bruhat}, which says that $\sigma\leq^*\tau$ is equivalent to $\sigma\leq_B\tau$. Finally, in Subsection \ref{sst:if_contained_then_star}, we prove that if $\Omega_\sigma\subseteq\overline{\Omega}_{\tau}$, then $\sigma\leq^*\tau$, see Proposition~\ref{prop:if_contained_then_star}. Thus, the conditions $\sigma\leq_B\tau$, $\Omega_\sigma\subseteq\overline{\Omega}_{\tau}$ and $\sigma\leq^*\tau$ are equivalent. This concludes the proof of our main result. Section~\ref{sect:remarks} contains some related facts and conjectures. In Subsection \ref{sst:dim_Omega}, we present a~formula for the dimension of the orbit $\Omega$ (see Theorem~\ref{theo:dim_Omega}). In Subsection \ref{sst:Schubert}, a conjectural approach to orbits associated with involutions in terms of tangent cones to Schubert varieties is described.

\medskip\textsc{Acknowledgements}. This work was done during my stay at Moscow State University. I~would like to express my gratitude to Professor E.B. Vinberg for his hospitality. Financial support from RFBR (grant no. 11--01--90703--mob${}_-$st) is gratefully acknowledged.

\sst\label{sst:A_n} Let $G=\GL_n(\Cp)$ be the general linear group, $B=B_n$ the subgroup of invertible upper-triangular matrices, and $U=U_n$ the unitriangular group, i.e., the unipotent radical of $B$. Let $\nt=\nt_n$ be the space of upper-triangular matrices with zeroes on the diagonal, and $\nt^*$ the dual space. Let $\Phi^+$ be the set of positive roots with respect to $B$. We identify $\Phi^+$ with the set $\{\epsi_i-\epsi_j,1\leq i<j\leq n\}$, where $\{\epsi_i\}_{i=1}^n$ is the standard basis of $\Rp^n$ (see, e.g., \cite{Bourbaki1}). Denote by $e_{i,j}$ the usual $(i,j)$th matrix unit, then $\{e_{\alpha},\alpha\in\Phi^+\}$ is a basis of $\nt^+$, where $e_{\epsi_i-\epsi_j}=e_{i,j}$. One can consider the dual basis $\{e_{\alpha}^*,\alpha\in\Phi^+\}$ of the dual space $\nt^*$.

The group $B$ acts on $\nt$ by conjugation, so one can consider the dual action of $B$ on $\nt^*$. By definition,
\begin{equation*}
\langle g.\lambda,x\rangle=\langle\lambda,g^{-1}xg\rangle,\text{ }g\in B,\text{ }x\in\nt,\text{ }\lambda\in\nt^*.
\end{equation*}
For a given $\lambda\in\nt^*$, let $\Omega_{\lambda}$ denote its orbit under this action. A subset $D\subset\Phi^+$ is called \emph{orthogonal} if it consists of pairwise orthogonal roots. To each orthogonal subset $D$ one can assign the element $$f_D=\sum_{\alpha\in D}e_{\alpha}^*\in\nt^*.$$ Let $W=W(A_{n-1})\cong S_n$ be the Weyl group of $G$. An involution $\sigma\in\Iu(A_{n-1})$ can be uniquely expressed as product of disjoint 2-cycles $\sigma=(i_1,j_1)\ldots(i_t,j_t)$, $i_l<j_l$ and $i_1<\ldots<i_t$. If we identify a~{trans\-po\-si\-tion} $(i,j)\in S_n$ with the reflection $r_{\epsi_i-\epsi_j}$ in the hyperplane orthogonal to the root $\epsi_i-\epsi_j$, then $$\sigma=\prod_{\alpha\in D}r_{\alpha},$$ where the orthogonal subset $D=\{\epsi_{i_1}-\epsi_{j_1},\ldots,\epsi_{i_t}-\epsi_{j_t}\}$ is called the \emph{support} of $\sigma$. We say that the $B$-orbit $\Omega_{\sigma}=\Omega_{f_D}$ is \emph{associated} with $\sigma$.

It is very convenient to identify $\nt^*$ with the space $\nt^t$ of lower-triangular matrices with zeroes on the diagonal by putting $e_{\alpha}^*=e_{\alpha}^t$. Under this identification, $$\langle\lambda,x\rangle=\tr\lambda x,\text{ }\lambda\in\nt^t,\text{ }x\in\nt.$$
For this reason, we will denote $\nt^t$ by $\nt^*$ and interpret it as the dual space of $\nt$. Note that if $g\in B$, $\lambda\in\nt^*$, then $$g.\lambda=(g\lambda g^{-1})_{\low},$$
where $A_{\low}$ denotes the strictly lower-triangular part of a matrix $A$. To each involution $\sigma\in\Iu(A_{n-1})$ one can assign the $0$--$1$ matrix $X_{\sigma}$ such that $(X_{\sigma})_{i,j}=1$ if and only if $\sigma(i)=j$.

\medskip\exam{It is convenient to draw a $0$--$1$ matrix $A$ as a \emph{rook placement} on the $n\times n$ board: by definition, there is a rook in the $(i,j)$th box if and only if $A_{i,j}=1$. For instance, let $n=6$, $\Supp{\sigma}=\{\epsi_1-\epsi_4,\epsi_3-\epsi_5\}$. On the picture below we draw $X_{\sigma}$ (rooks are marked by $\otimes$'s).
\begin{equation*}\predisplaypenalty=0X_{\sigma}=
\mymatrix{
\pho& \pho& \pho& \otimes& \pho& \pho\\
\pho& \otimes& \pho& \pho& \pho& \pho\\
\pho& \pho& \pho& \pho& \otimes& \pho\\
\otimes& \pho& \pho& \pho& \pho& \pho\\
\pho& \pho& \otimes& \pho& \pho& \pho\\
\pho& \pho& \pho& \pho& \pho& \otimes\\
}\
\end{equation*}}

\medskip To each $0$--$1$ matrix $A$ one can assign the matrix $R(A)$ by putting
$$R(A)_{i,j}=\rk\pi_{i,j}(A),$$
where $\pi_{i,j}$ denotes the lower-left $(n-i+1)\times j$ submatrix of $A$. (In other words, $R(A)_{i,j}$ is just the number of rooks located
non-strictly to the South-West of the $(i,j)$th box.) In particular, we set $R_{\sigma}=R(X_{\sigma})$ and $R_{\sigma}^*=(R_{\sigma})_{\low}$, the lower-triangular part of $R_{\sigma}$. Suppose $A$ and $B$ are matrices with integer entries. We write $A\leq B$ if $A_{i,j}\leq B_{i,j}$ for all $i,j$.
Denote by $\leq_B$ the usual Bruhat--Chevalley order on $S_n$. Denote also by $\overline{Z}$ the closure of a subset $Z\subseteq\nt^*$ with respect to Zariski topology. We have the following description of the incidences among the closures of $B$-orbits associated with involutions.

\medskip\theop{Let $\sigma,\tau\in\Iu(A_{n-1})$. The following conditions are equivalent:
\begin{equation*}
\begin{split}
&\text{\rm{i) }}\sigma\leq_B\tau;\\
&\text{\rm{ii) }}R_{\sigma}\leq R_{\tau};\\
&\text{\rm{iii) }}R_{\sigma}^*\leq R_{\tau}^*;\\
&\text{\rm{iv) }}\Omega_{\sigma}\subseteq\overline{\Omega}_{\tau}.\\
\end{split}
\end{equation*}\label{mtheo_A_n}}{i) $\Leftrightarrow$ ii). See, e.g., \cite[Theorem 1.6.4]{Incitti1}.

ii) $\Rightarrow$ iii) is evident. The converse follows from \cite[Theorem 1.10]{Ignatev3}.

iii) $\Leftrightarrow$ iv). See \cite[Theorem 1.7]{Ignatev3}.}

\sst\label{sst:C_n} From now on, let $G=\Sp_{2n}(\Cp)$ be the \emph{symplectic} group, i.e.,
\begin{equation*}
\begin{split}
&G=\{X\in\GL_{2n}(\Cp)\mid X^tJX=J\},\text{ where }\\
&J=\begin{pmatrix}0&s\\-s&0\end{pmatrix}.
\end{split}
\end{equation*}
Here $s$ denotes $n\times n$ matrix with $1$'s on the antidiagonal and zeroes elsewhere. Let $B=B_n\cap G$ be the Borel subgroup of $G$ consisting of all upper-triangular matrices from $G$. The unipotent radical of $B$ coincides with the group $U=U_n\cap G$ of all upper-triangular matrices from $G$ with $1$'s on the diagonal.

The Lie algebra $\gt$ of $G$ is the \emph{symplectic} algebra $\spt_{2n}(\Cp)=\{X\in\Mat_{2n}(\Cp)\mid X^tJ+JX=0\}$. The Lie algebra of $B$ coincides with the space $\bt$ of all upper-triangular matrices from $\gt$. The space $\nt$~of matrices from $\bt$ with zeroes on the diagonal is the Lie algebra of $U$. Let $\Phi^+$ be set of positive roots with respect to $B$. We identify $\Phi^+$ with the set $$\{\epsi_i\pm\epsi_j,1\leq i<j\leq n\}\cup\{2\epsi_i,1\leq i\leq n\},$$ where $\{\epsi_i\}_{i=1}^n$ is the standard basis of $\Rp^n$ (see, e.g., \cite{Bourbaki1}).

For a given $\alpha\in\Phi^+$, put
\begin{equation*}
e_{\alpha}=\begin{cases}e_{i,j}-e_{-j,-i},\text{ if }\alpha=\epsi_i-\epsi_j,\\
e_{i,-j}+e_{j,-i},\text{ if }\alpha=\epsi_i+\epsi_j,\\
e_{i,-i},\text{ if }\alpha=2\epsi_i.\\
\end{cases}
\end{equation*}
Here we index rows and columns of any $2n\times2n$ matrix by the numbers $1,2,\ldots,n,-n,-n+1,\ldots,-1$ and denote by $e_{i,j}$ the $(i,j)$th matrix unit. The set $\{e_{\alpha},\alpha\in\Phi^+\}$ is a basis of $\nt$, so one can consider the dual basis $\{e_{\alpha}^*,\alpha\in\Phi^+\}$ of the dual space $\nt^*$.

Since $\nt$ is invariant under the adjoint action of $B$ on $\bt$, one can consider the dual action of $B$ on~$\nt^*$. By definition, if $g\in B$, $\lambda\in\nt^*$ and $x\in\nt$, then $$\langle g.\lambda,x\rangle=\langle\lambda,\Ad_g^{-1}(x)\rangle,$$ where $\Ad$ denotes the adjoint action. (In fact, $\Ad_g(x)=gxg^{-1}$.) For a given $\lambda\in\nt^*$, by $\Omega_{\lambda}$ (resp. by~$\Theta_{\lambda}$) we denote the $B$-orbit (resp. the $U$-orbit) of $\lambda$.

A subset $D\subset\Phi^+$ is called \emph{orthogonal} if it consists of pairwise orthogonal roots. To each orthogonal subset $D\subset\Phi$ and each map $\xi\colon D\to\Cp^{\times}$ one can assign the elements of $\nt^*$ of the form $$f_D=\sum_{\alpha\in D}e_{\alpha}^*,\text{ }f_{D,\xi}=\sum_{\alpha\in D}\xi(\alpha)e_{\alpha}^*.$$ (If $D=\emptyset$, then $f_{D,\xi}=0$.) Evidently, $f_D=f_{D,\xi_0}$, where $\xi_0$ sends all roots from $D$ to $1$. Put $\Omega_D=\Omega_{f_D}$, $\Theta_D=\Theta_{f_D}$ and $\Theta_{D,\xi}=\Theta_{f_{D,\xi}}$. Clearly, $\Theta_D\subseteq\Omega_D$; Lemma \ref{lemm:Omega_ravno_union_Theta} shows that in fact $\Omega_D=\bigcup\Theta_{D,\xi}$, where the union is taken over all maps $\xi\colon D\to\Cp^{\times}$.

\medskip\nota{Orbits associated with orthogonal subsets play an important role in representation theory of $U$. They were studied by the author in \cite{Ignatev1}. (See also \cite{Ignatev2} for further examples and generalizations to other unipotent algebraic groups.)}

\medskip Let $\sigma\in\Iu(\Phi)$ be an involution from $W$, the Weyl group of $\Phi$. An orthogonal subset $D\subset\Phi^+$ is called a \emph{support} of $\sigma$, if $\sigma=\prod_{\alpha\in D}r_{\alpha}$ and there are no $\alpha,\beta\in D$ such that $\alpha-\beta\in\Phi^+$. Here $r_{\alpha}\in W$ denotes the reflection in the hyperplane orthogonal to a given root $\alpha$, and the product is taken in any fixed order. One can easily check that there exists exactly one support of $\sigma$ among all orthogonal subsets of $\Phi^+$. We denote it by $\Supp{\sigma}$. We put also $\Omega_{\sigma}=\Omega_{\Supp{\sigma}}$ and $f_{\sigma}=f_{\Supp{\sigma}}$.

\medskip\exam{Let $D=\{\epsi_1-\epsi_2,\epsi_1+\epsi_2\}$, $D'=\{2\epsi_1,2\epsi_2\}$. Then $\sigma=\prod_{\alpha\in D}r_{\alpha}=\prod_{\alpha\in D'}r_{\alpha}$, but $(\epsi_1+\epsi_2)-(\epsi_1-\epsi_2)=2\epsi_2\in\Phi^+$, so $\Supp{\sigma}=D'$, not $D$.}

\medskip\defi{We say that the $B$-orbit $\Omega_{\sigma}$ is \emph{associated} with the involution $\sigma$.\label{defi:B_orbit_assoc_invol}}

\medskip\nota{i) If $D$ is an orthogonal subset of $C_n^+$ such that $\alpha-\beta\in\Phi^+$  for some $\alpha,\beta\in D$, then, by \cite[Proposition 2.1]{Ignatev1}, $\Theta_{D,\xi}=\Theta_{D_1,\xi_1}$ for any $\xi\colon D\to\Cp^{\times}$, where $D_1=D\setminus\{\beta\}$ and $\xi_1=\xi\mathbin{\mid}_{D_1}$. Applying Lemma \ref{lemm:Omega_ravno_union_Theta}, we see that $\Omega_D=\Omega_{D_1}$, while $\sigma\neq\sigma_1$, where $\sigma=\prod_{\alpha\in D}r_{\alpha}$, $\sigma_1=\prod_{\alpha\in D_1}r_{\alpha}$. Note, however, that if $\alpha-\beta\in\Phi^+$, then $\alpha=\epsi_i+\epsi_j$, $\beta=\epsi_i-\epsi_j$ for some $i,j$. Thus, $\Supp{\sigma}=D'$, where $D'$ is obtained from $D$ by replacing each pair $\{\epsi_i-\epsi_j,\epsi_i+\epsi_j\}$ by $\{2\epsi_i,2\epsi_j\}$.

ii) In fact, we do \emph{not} know how to define the support of an involution for other root systems. For instance, let $\Phi=D_n$, $n\geq4$, then $\Phi^+$ can be identified with $\{\epsi_i\pm\epsi_j,1\leq i<j\leq n\}$. Put
\begin{equation*}
\begin{split}
D&=\{\epsi_1-\epsi_2,\epsi_1+\epsi_2,\epsi_3-\epsi_4,\epsi_3+\epsi_4\},\\
D'&=\{\epsi_1-\epsi_3,\epsi_1+\epsi_3,\epsi_2-\epsi_4,\epsi_2+\epsi_4\}.
\end{split}
\end{equation*}
Then $\sigma=\prod_{\alpha\in D}r_{\alpha}=\prod_{\alpha\in D'}r_{\alpha}\in\Iu(D_n)$, so we have two candidates for the role of $\Supp{\sigma}$ and we do not know how to choose one of them.}

\medskip It is very convenient to identify $\nt^*$ with the space $\nt^t$ of lower-triangular matrices from $\gt$ by putting $e_{\alpha}^*=e_{\alpha}^t$. Under this identification, $$\langle\lambda,x\rangle=\tr\lambda'x,\text{ }\lambda\in\nt^t,\text{ }x\in\nt.$$
Here we put $\Phi_0^+=\{\epsi_i\pm\epsi_j,1\leq i<j\leq n\}$, $\Phi_1^+=\{2\epsi_i,1\leq i\leq n\}$; evidently, $\Phi^+=\Phi_0^+\cup\Phi_1^+$, so $\nt=\nt_0\oplus\nt_1$ as vector spaces, where $\nt_0$ (resp. $\nt_1$) is spanned by $e_{\alpha}$, $\alpha\in\Phi_0^+$ (resp. $\alpha\in\Phi_1^+$); if $\lambda=\lambda_0+\lambda_1$, $\lambda_0\in\nt_0$, $\lambda_1\in\nt_1$, then, by definition, $\lambda'=1/2\lambda_0+\lambda_1$. For this reason, we will denote $\nt^t$ by $\nt^*$ and interpret it as the dual space of $\nt$. Note that if $g\in B$, $\lambda\in\nt^*$, then $$g.\lambda=(g\lambda g^{-1})_{\low},$$
where $A_{\low}$ denotes the strictly lower-triangular part of a matrix $A$ (see, e.g., \cite[p. 410]{AndreNeto1}).

\medskip\lemmp{Let $D$ be an orthogonal subset of $C_n^+$. Then\footnote{Cf. \cite[Lemma 2.1]{Ignatev3}} $\Omega_D=\bigcup\Theta_{D,\xi}$, where the union is taken over all maps $\xi\colon D\to\Cp^{\times}$.\label{lemm:Omega_ravno_union_Theta}}{It is well-known (see, f.e., \cite[Subsection 15.1]{Humphreys1}) that the \emph{exponential} map $$\exp\colon\nt\to U\colon x\mapsto\sum_{i=0}^{\infty}\dfrac{x^i}{i!}$$ is well-defined; in fact, it is an isomorphism of affine varieties. For a given $\alpha\in\Phi^+$, $s\in\Cp^{\times}$, put
\begin{equation*}\predisplaypenalty=0
\begin{split}
&x_{\alpha}(s)=\exp(se_{\alpha})=1+se_{\alpha},\text{ }x_{-\alpha}(s)=x_{\alpha}(s)^t,\\
&w_{\alpha}(s)=x_{\alpha}(s)x_{-\alpha}(-s^{-1})x_{\alpha}(s),\text{ }h_{\alpha}(s)=w_{\alpha}(s)w_{\alpha}(1)^{-1}.
\end{split}
\end{equation*}
Then $h_{\alpha}(s)$ is a diagonal matrix from $B$.

Let $\xi\colon D\to\Cp^{\times}$ be a map. Suppose $\alpha\in D$. Pick a number $s\in\Cp^{\times}$ and put
\begin{equation*}
s'=\begin{cases}s^{-1},&\text{if }\alpha=\epsi_i+\epsi_j,\\
s,&\text{if }\alpha=\epsi_i-\epsi_j,\\
\sqrt{s^{-1}},&\text{if }\alpha=2\epsi_i,\\
\end{cases}
\end{equation*}
and $\alpha'=2\epsi_i$ in all cases above. (Here $\sqrt{s^{-1}}$ is a number such that $(\sqrt{s^{-1}})^2=s^{-1}$.) One can easily check by straightforward matrix calculations that $$h_{\alpha'}(s').f_{D,\xi}=\sum_{\beta\in D,\beta\neq\alpha}\xi(\beta)e_{\beta}^*+s\xi(\alpha)e_{\alpha}^*.$$ Hence $$\left(\prod_{\alpha\in D}h_{\alpha'}(\xi(\alpha)')\right).f_D=f_{D,\xi},$$
so $\Theta_{D,\xi}\subseteq\Omega_{D}.$

On the other hand, let $h\in H$, where $H$ is the group of diagonal matrices from $G$. We claim that $h.f_{D,\xi}=f_{D,\xi'}$ for some $\xi'$. Indeed, since $H$ is generated by $h_{\alpha}(s)$'s, $\alpha\in\Phi^+$, $s\in\Cp^{\times}$, we can assume without loss of generality that $h=h_{\alpha}(s)$ for some $\alpha$ and $s$. But in this case the statement follows immediately from the above.

The group $B$ is isomorphic as an algebraic group to the semi-direct product $U\rtimes H$ \cite[Subsection 19.1]{Humphreys1}. In particular, for a~given $g\in B$, there exist unique $u\in U$, $h\in H$ such that $g=uh$. If $\xi\colon D\to\Cp^{\times}$ is the map such that $h.f_D=f_{D,\xi}$, then $g.f_D=u.f_{D,\xi}\in\Theta_{D,\xi}$. This concludes the proof.}

\sect{Proof of the main theorem}\label{sect:proof_mtheo}
\sst\label{sst:if_Bruhat_then_contained} In this subsection, we will prove that if $\sigma,\tau\in\Iu(C_n)$ and $\sigma$ is less or equal to $\tau$ with respect to the Bruhat--Chevalley order, then $\Omega_{\sigma}$ is contained in $\overline{\Omega}_{\tau}$, the Zariski closure of $\Omega_{\tau}$. We will denote the set of \emph{fundamental} roots $\{\epsi_1-\epsi_2,\ldots,\epsi_{n-1}-\epsi_n,2\epsi_n\}$ by $\Pi$. A reflection $r_{\alpha}$ is called \emph{fundamental} if $\alpha\in\Pi$. An expression of $w\in W$ as a product of fundamental reflections is called \emph{reduced} if it has the minimal length among all such expressions. The length $l(w)$ of a reduced expression is called the \emph{length} of $w$. Let $w=r_{\alpha_1}\ldots r_{\alpha_l}$ be a reduced expression. By definition of the \emph{Bruhat--Chevalley order}~$\leq_B$ on $W$, $$\{w'\in W\mid w'\leq_B w\}=\{r_{\alpha_1}\ldots r_{\alpha_{i_t}},\text{ }t\leq l,\text{ }i_1<\ldots<i_t\}.$$

Let $\sigma,\tau\in\Iu(\Phi)$. We say that $\tau$ \emph{covers} $\sigma$ and write $\sigma\tri\tau$ if $\tau<_B\sigma$ and there are no $w\in\Iu(\Phi)$ such that $\tau<_Bw<_B\sigma$. In \cite{Incitti1}, F. Incitti studied the restriction of $\leq_B$ to $\Iu(\Phi)$ from a combinatorial point of view. In particular, for a given involution $\tau$, he described the set $L(\tau)=\{\sigma\in\Iu(\Phi)\mid\sigma\tri\tau\}$ (see \cite[pp. 75--81]{Incitti1}). We reformulate his description in our terms in Appendix.

\medskip\propp{Let $\sigma,\tau$ be involutions in $W(C_n)$. If $\sigma\leq_B\tau$, then $\Omega_{\sigma}\subseteq\overline{\Omega}_{\tau}$.\label{prop:if_Bruhat_then_contained}}{There exist $\tau_1,\ldots,\tau_r\in\Iu(\Phi)$ such that $\tau_1=\sigma$, $\tau_r=\tau$ and $\tau_i\tri\tau_{i+1}$ for all $1\leq i<r$, so we can assume without loss of generality that $\sigma\tri\tau$. It is well-known that $\overline{\Omega}_{\tau}$ coincides with $\overline{\Omega}_{\tau}^{\Cp}$, the closure of $\Omega_{\tau}$ with respect to the complex topology, so it suffice to construct $g(s)\in B$, $s\in\Cp^{\times}$, such that $g(s).f_{\sigma}\to f_{\tau}$ as $s\to0$.


The proof is case-by-case. For example, let
\begin{equation*}
\begin{split}
&\Supp{\sigma}\setminus\Supp{\tau}=\{\epsi_i+\epsi_j,2\epsi_k\},\\
&\Supp{\tau}\setminus\Supp{\sigma}=\{2\epsi_i,\epsi_k+\epsi_j\}
\end{split}
\end{equation*}
for some $1\leq i<k<j\leq n$ (see Case 5 in Appendix). Put
$$g(s)=h_{\epsi_i-\epsi_k}(s^{-1})\cdot x_{\epsi_i-\epsi_k}(s)\cdot x_{\epsi_k-\epsi_j}((2s^2)^{-1})\cdot x_{\epsi_i-\epsi_j}(-s^{-1})$$
and $f=g(s).f_{\tau}$. One can easily check by straightforward matrix calculations that
\begin{equation*}
f(e_{\alpha})=\begin{cases}1,&\text{if either }\alpha=\epsi_i+\epsi_j\text{ or }\alpha=2\epsi_k,\\
0,&\text{if }\alpha=\epsi_k+\epsi_j,\\
-s,&\text{if }\alpha=\epsi_i+\epsi_k,\\
s^2,&\text{if }\alpha=2\epsi_i,\\
f_{\tau}(e_{\alpha})&\text{otherwise}.\\
\end{cases}
\end{equation*}
Thus, $f\to f_{\sigma}$ as $s\to0$.

All other cases can be considered similarly, see Appendix.}

\sst\label{sst:star_equiv_Bruhat} Let $S_{\pm n}$ be the symmetric group on $2n$ letters $\{1,2,\ldots,n,-n,-n+1,\ldots,-1\}$. Put
$$W'=\{w\in S_{\pm n}\mid w(-i)=-w(i)\text{ for all }1\leq i\leq n\}.$$
It is known that the map $W\to W'\colon w\mapsto w'$ defined by
\begin{equation*}
\begin{split}
&r_{\epsi_i-\epsi_j}'=(i,j)(-i,-j),\\
&r_{\epsi_i+\epsi_j}'=(i,-j)(-i,j),\\
&r_{2\epsi_i}'=(i,-i)
\end{split}
\end{equation*}
is an isomorphism of groups. Note that if $\sigma\in\Iu(C_n)$, then $\sigma'$ is an involution in $S_{\pm n}$.

To each involution $\sigma\in\Iu(C_n)$ one can assign the $0$--$1$ matrix $X_{\sigma}$ such that $(X_{\sigma})_{i,j}=1$ if and only if $\sigma'(i)=j$ (recall that we index rows and columns of a $2n\times 2n$ matrix by the numbers $1,\ldots,n,-n,\ldots,1$).

\medskip\exam{It is convenient to draw a $0$--$1$ matrix $A$ as a rook placement on the $2n\times 2n$ board: by definition, there is a rook in the $(i,j)$th box if and only if $A_{i,j}=1$. For instance, let $n=4$, $\Supp{\sigma}=\{\epsi_1-\epsi_4,2\epsi_2\}$. On the picture below we draw $X_{\sigma}$ (rooks are marked by $\otimes$'s).
\begin{equation*}\predisplaypenalty=0X_{\sigma}=
\mymatrixII{ \pho& \pho& \pho& \otimes& \pho& \pho& \pho& \pho\\
\pho& \pho& \pho& \pho& \pho& \pho& \otimes& \pho\\
\pho& \pho& \otimes& \pho& \pho& \pho& \pho& \pho\\
\otimes& \pho& \pho& \pho& \pho& \pho& \pho& \pho\\
\pho& \pho& \pho& \pho& \pho& \pho& \pho& \otimes\\
\pho& \pho& \pho& \pho& \pho& \otimes& \pho& \pho\\
\pho& \otimes& \pho& \pho& \pho& \pho& \pho& \pho\\
\pho& \pho& \pho& \pho& \otimes& \pho& \pho& \pho\\
}\
\end{equation*}}

\medskip As above, to each $0$--$1$ matrix $A$ we assign the matrix $R(A)$ by putting
$$R(A)_{i,j}=\rk\pi_{i,j}(A),$$
where $\pi_{i,j}$ sends $A$ to its submatrix with rows $i,\ldots,-1$ and columns $1,\ldots,j$. (In other words, $R(A)_{i,j}$ is just the number of rooks located
non-strictly to the South-West of the $(i,j)$th box.) In particular, we set $R_{\sigma}=R(X_{\sigma})$ and $R_{\sigma}^*=(R_{\sigma})_{\low}$, the lower-triangular part of $R_{\sigma}$.


Suppose $A$ and $B$ are matrices with integer entries. As above, we write $A\leq B$ if $A_{i,j}\leq B_{i,j}$ for all $i,j$. Let $\sigma,\tau\in\Iu(C_n)$. By \cite[Theorem 1.6.7]{Incitti1}, $$\sigma\leq_B\tau\text{ if and only if }R_{\sigma}\leq R_{\tau}.$$
(Note that Incitti use another order of fundamental roots.) Let us define another partial order on~$\Iu(C_n)$. Namely, we put
$$\sigma\leq^*\tau\text{ if }R_{\sigma}^*\leq R_{\tau}^*.$$
Clearly, $\sigma\leq_B\tau$ implies $\sigma\leq^*\tau$. But the converse follows immediately from Theorem~\ref{mtheo_A_n}, since $\sigma'$~and~$\tau'$~are involutions in $S_{\pm n}$. This proves

\medskip\post{Let $\sigma,\tau$ be involutions in $W(C_n)$. Then $\sigma\leq_B\tau$ if and only if $\sigma\leq^*\tau$.\label{prop:star_equiv_Bruhat}}

\sst\label{sst:if_contained_then_star} Suppose $\sigma,\tau\in\Iu(C_n)$. To conclude the proof of Theorem \ref{mtheo}, it remains to check that if $\Omega_{\sigma}\subseteq\overline{\Omega}_{\tau}$, then $\sigma\leq^*\tau$. To do this, we need the following

\medskip\lemmp{Let $\sigma\in\Iu(C_n)$. Then\footnote{Cf. \cite[Lemma 2.2]{Ignatev3}} $\rk\pi_{i,j}(\lambda)=(R_{\sigma}^*)_{i,j}$ for all $\lambda\in\Omega_{\sigma}$.\label{lemm:rank_through_orbit}}{Fist, note that $(R_{\sigma}^*)_{i,j}=\rk\pi_{i,j}(f_\sigma)$. By Lemma \ref{lemm:Omega_ravno_union_Theta}, $\Omega_{\sigma}=\bigcup_{\xi\colon D\to\Cp^{\times}}\Omega_{D,\xi}$, where ${D=\Supp{\sigma}}$. Let $\xi\colon D\to\Cp^{\times}$ be a map. Since $\rk\pi_{i,j}(f_{D,\xi})=\rk\pi_{i,j}(f_{\sigma})=(R_{\sigma}^*)_{i,j}$ for all $i,j$, it suffice to check that $\rk\pi_{i,j}(\lambda)=\rk\pi_{i,j}(u.\lambda)$ for all $u\in U$, $\lambda\in\nt^*$. But this follows immediately from the proof of \cite[Lemma 2.2]{Ignatev3}, because $u$ is an upper-triangular matrix with $1$'s on the diagonal, $\lambda$ is a~lower-triangular matrix with zeroes on the diagonal, and $u.\lambda=(u\lambda u^{-1})_{\low}$.}

\medskip Things now are ready to prove

\medskip\propp{Let $\sigma,\tau$ be involutions in $W(C_n)$. If $\Omega_{\sigma}\subseteq\overline{\Omega}_{\tau}$, then\footnote{Cf. \cite[Proposition 2.3]{Ignatev3}} $\sigma\leq^*\tau$\label{prop:if_contained_then_star}}{Suppose
$\sigma\nleq^*\tau$. This means that there exist $i,j$ such that $(R_{\sigma}^*)_{i,j}>(R_{\tau}^*)_{i,j}$.
Denote
\begin{equation*}
Z=\{f\in\nt^*\mid\rk\pi_{r,s}(f)\leq(R_{\tau}^*)_{r,s}\text{ for
all }r,s\}.
\end{equation*}
Clearly, $Z$ is closed with respect to the Zariski topology. Lemma
\ref{lemm:rank_through_orbit} shows that $\Omega_{\tau}\subseteq
Z$, so $\overline{\Omega}_{\tau}\subseteq Z$. But
$f_{\sigma}\notin Z$, hence $\Omega_{\sigma}\nsubseteq Z$, a
contradiction.}

\medskip The proof of Theorem \ref{mtheo} is complete.

\sect{Concluding remarks}\label{sect:remarks}

\sst\label{sst:dim_Omega} Let $\sigma\in\Iu(C_n)$. Being an orbit of a connected unipotent group on an affine variety $\nt^*$, $\Omega_{\sigma}$ is a closed subvariety of $\nt^*$. In this subsection, we present a formula for the dimension of $\Omega_{\sigma}$. Recall the definition of the length $l(w)$ of an element $w\in W$.

\medskip\theop{Let $\sigma\in W(C_n)$ be an involution. Then\footnote{Cf. \cite[Proposition 4.1]{Ignatev3}} $\dim\Omega_{\sigma}=l(\sigma)$.\label{theo:dim_Omega}}{Put $D=\Supp{\sigma}$. We claim that if $\xi_1,\xi_2$ are two distinct maps from $D$ to $\Cp^{\times}$, then $\Theta_{D,\xi_1}\neq\Theta_{D,\xi_2}$. Indeed, let $\wt U=U_n$ be the unitriangular group (i.e., the group of all upper-triangular matrices with $1$'s on the diagonal). Since $\sigma'$ is an involution in $S_{\pm n}$, \cite[Theorem 1.4]{Panov} implies that $\wt\Theta_{D,\xi_1}\neq\wt\Theta_{D,\xi_2}$, where $\wt\Theta_{D,\xi_1}$ (resp.\ $\wt\Theta_{D,\xi_2}$) denotes the $\wt U$-orbit of $f_{D,\xi_1}$ (resp.\ of $f_{D,\xi_2}$) under the action of $\wt U$ on the space of all lower-triangular matrices with zeroes on the diagonal defined by the formula $$u.\lambda=(u\lambda u^{-1})_{\low},\text{ }u\in\wt U,\text{ }\lambda\in\wt\nt^*.$$
Since $U\subseteq\wt U$, one has $\Theta_{D,\xi_1}\subseteq\wt\Theta_{D,\xi_1}$ and $\Theta_{D,\xi_2}\subseteq\wt\Theta_{D,\xi_2}$, hence $\Theta_{D,\xi_1}\neq\Theta_{D,\xi_2}$, as required.

Let $Z_B=\Stab_Bf_{\sigma}$ be the stabilizer of $f_{\sigma}$ in $B$. One has $$\dim\Omega_{\sigma}=\dim B-\dim Z_B.$$
Recall that $B\cong U\rtimes H$ as algebraic groups. It was shown in the proof of Lemma \ref{lemm:Omega_ravno_union_Theta} that if $h\in H$, then there exists $\xi\colon D\to\Cp^{\times}$ such that $h.f_{\sigma}=f_{D,\xi}$. Hence if $g=uh\in Z_B$, then $$f_{\sigma}=(uh).f_{\sigma}=u.f_{D,\xi},$$
so $f_{\sigma}\in\Theta_{D,\xi}$. In follows from the first paragraph of the proof that $f_{\sigma}=f_{D,\xi}$. This means that the map $$Z_U\times Z_H\to Z_B\colon (u,h)\mapsto uh$$ is an isomorphism of algebraic varieties, where $Z_U=\Stab_Uf_{\sigma}$ (resp. $Z_H=\Stab_Hf_{\sigma}$) is the stabilizer of $f_{\sigma}$ in $U$ (resp. in $H$). Hence $$\dim Z_B=\dim Z_U+\dim Z_H.$$

By \cite[Theorem 1.2]{Ignatev1}, $\dim\Theta_{D}=l(\sigma)-|D|$, so $$\dim Z_U=\dim U-\dim\Theta_D=\dim U-l(\sigma)+|D|.$$ On the other hand, put $X=\bigcup_{\xi\colon D\to\Cp^{\times}}\{f_{D,\xi}\}$. In follows from Lemma \ref{lemm:Omega_ravno_union_Theta} and the first paragraph of the proof that $X=\{h.f_{\sigma},h\in H\}$, the $H$-orbit of $f_{\sigma}$. Consequently, $$\dim Z_H=\dim H-\dim X=\dim H-|D|,$$ because $X$ is isomorphic as affine variety to the product of $|D|$ copies of $\Cp^{\times}$. Thus,
\begin{equation*}
\begin{split}
\dim\Omega_{\sigma}&=\dim B-\dim Z_B=(\dim U+\dim H)-(\dim Z_U+\dim Z_H)\\
&=\dim U+\dim H-(\dim U-l(\sigma)-|D|)-(\dim H-|D|)=l(\sigma).
\end{split}
\end{equation*} The proof is complete.}

\sst\label{sst:Schubert} In the remainder of the paper, we briefly discuss a conjectural geometrical approach to orbits associated with involutions in terms of tangent cones to Schubert varieties. Recall that $W$ is isomorphic to $N_G(H)/H$, where $N_G(H)$ is the normalizer of $H$ in $G$. The \emph{flag variety} $\Fo=G/B$ can be decomposed into the union $\Fo=\bigcup_{w\in W}X_w^{\circ}$, where $X_w^{\circ}=B\dot wB/B$
is called the \emph{Schubert cell}. (Here $\dot w$ is a representative of $w$ in $N_G(H)$.) By definition, the \emph{Schubert variety} $X_w$ is the closure of
$X_w^{\circ}$ in~$\Fo$ with respect to Zariski topology. Note that $p=X_{\id}=B/B$ is contained in $X_w$ for all $w\in W$. One has $X_w\subseteq X_{w'}$ if and only if $w\leq_Bw'$. Let $T_w$ be the tangent space and $C_w$ the tangent cone to~$X_w$ at
the point $p$ (see \cite{BileyLakshmibai1} for detailed constructions); by definition, $C_w\subseteq T_w$,
and if $p$ is a~regular point of~$X_w$, then $C_w=T_w$. Of course, if $w\leq_Bw'$, then $C_w\subseteq C_{w'}$.

Let $T=T_p\Fo$ be the tangent space to $\Fo$ at $p$. It can be naturally identified with $\nt^*$ by the following way: since $\Fo=G/B$, $T$ is isomorphic to the factor $\gt/\bt\cong\nt^*$. Next, $B$ acts on $\Fo$ by conjugation. Since $p$ is invariant under this action, the action on $T=\nt^*$ is induced. One can check that this action coincides with the action of $B$ on $\nt^*$ defined above. The tangent cone $C_w\subseteq T_w\subseteq T=\nt^*$
is $B$-invariant, so it splits into a union of $B$-orbits. Furthermore, $\overline{\Omega}_{\sigma}\subseteq C_{\sigma}$ for all $\sigma\in\Iu(C_n)$.

It is well-known that $C_w$ is a subvariety of $T_w$ of dimension $\dim C_w=l(w)$ \cite[Chapter 2, Section~2.6]{BileyLakshmibai1}. Let $\sigma\in\Iu(C_n)$. Since $\Omega_{\sigma}$ is irreducible, $\overline{\Omega}_{\sigma}$ is. Theorem \ref{theo:dim_Omega} implies $\dim\overline{\Omega}_{\sigma}=\dim\Omega_{\sigma}=l(\sigma)$, so $\overline{\Omega}_{\sigma}$ is an irreducible component of $C_{\sigma}$ of maximal dimension.

\medskip\hypo{Let $\sigma\in W(C_n)$ be an involution. Then\footnote{Cf. \cite[Conjecture 1.11]{Ignatev3}.} the closure of the $B$-orbit $\Omega_{\sigma}$ coincides with the tangent cone $C_w$ to the Schubert variety $X_w$ at the point $p=B/B$.}

\medskip Note that this conjecture implies that if $\sigma\leq_B\tau$, then $\Omega_{\sigma}\subseteq\overline{\Omega}_{\tau}$.

\bigskip\begin{center}{\large\textbf{Appendix}}
\end{center}

Let $\tau\in\Iu(C_n)$, $L(\tau)=\{\sigma\in\Iu(C_n)\mid\sigma\tri\tau\}$. Here we reformulate Incitti's description of $L(\tau)$ used in the proof of Proposition \ref{prop:if_Bruhat_then_contained}. Precisely, for each pair $(\sigma,\tau)$, $\sigma\in L(\tau)$, we construct the matrix $g(s)\in B$, $s\in\Cp^{\times}$, such that $f=g(s).f_{\tau}\to f_{\sigma}$ as $s\to0$. By \cite[p. 96]{Incitti1}, $\sigma\in L(\tau)$ if and only of $(\sigma,\tau)$ is one of the pairs described below. This concludes the proof of Proposition \ref{prop:if_Bruhat_then_contained}. In the second column of the table below, we indicate the number of a case referred to Incitti's paper \cite{Incitti1}. (Note, however, that Incitti use another order of fundamental roots.) Note that there are some additional conditions to be satisfied by $(\sigma,\tau)$, $\sigma\in L(\tau)$ (see \cite[pp. 76--81]{Incitti1}), but $g(s)$ does not depend on them. In particular, all $\epsi_i$ occurring in the third or in the fourth column of the table must satisfy $(\beta,\epsi_i)=0$ for all $\beta\in D_{\sigma}\cap D_{\tau}$. We put here $D_{\sigma}=\Supp{\sigma}$ and $D_{\tau}=\Supp{\tau}$. We denote by $I$ a fixed complex number such that $I^2=-1$. For simplicity, if $f(e_{\alpha})=f_{\tau}(e_{\alpha})$, then we omit $\alpha$ in the last column of the table.

\begin{longtable}{||l|l|l|l|l|l||}\hline
&&$D_{\sigma}\setminus D_{\tau}$&$D_{\tau}\setminus D_{\sigma}$&$g(s)$&$f(e_{\alpha})$, $\alpha\in\Phi^+$\\
\hline\hline
1\rule{0pt}{12pt}&\small{$(A1,M1)$}&
\begin{tabular}[t]{l}$\emptyset$\end{tabular}&
\begin{tabular}[t]{l}$2\epsi_i$\end{tabular}&
\begin{tabular}[t]{l}$h_{2\epsi_i}(s^{-1})$\end{tabular}&
\begin{tabular}[t]{ll}$s^2$,&if $\alpha=2\epsi_i$.
\end{tabular}\\
\hline
2\rule{0pt}{12pt}&\small{$(A1,M5)$}&
\begin{tabular}[t]{l}$\epsi_i+\epsi_j$,\\
$i<j$
\end{tabular}&
\begin{tabular}[t]{l}$2\epsi_i$, $2\epsi_j$\end{tabular}&
\begin{tabular}[t]{l}$h_{2\epsi_i}(s^{-1})\times$\\
$h_{2\epsi_j}(-Is)\times$\\
$x_{\epsi_i-\epsi_j}(I)$
\end{tabular}&
\begin{tabular}[t]{ll}$s^2$,&if $\alpha=2\epsi_i$,\\
$0$,&if $\alpha=2\epsi_j$,\\
$1$,&if $\alpha=\epsi_i+\epsi_j$.
\end{tabular}\\
\hline
3\rule{0pt}{12pt}&\small{$(A1,M6)$}&
\begin{tabular}[t]{l}$\epsi_i-\epsi_j$,\\
$i<j$
\end{tabular}&
\begin{tabular}[t]{l}$2\epsi_i$\end{tabular}&
\begin{tabular}[t]{l}$h_{2\epsi_i}(s^{-1})\times$\\
$x_{\epsi_i+\epsi_j}(s^{-1})$
\end{tabular}&
\begin{tabular}[t]{ll}$s^2$,&if $\alpha=2\epsi_i$,\\
$1$,&if $\alpha=\epsi_i-\epsi_j$.
\end{tabular}\\
\hline
4\rule{0pt}{12pt}&\small{$(A2,M2)$}&
\begin{tabular}[t]{l}$2\epsi_j$
\end{tabular}&
\begin{tabular}[t]{l}$2\epsi_i$,\\
$i<j$\end{tabular}&
\begin{tabular}[t]{l}$h_{2\epsi_i}(s^{-1})\times$\\
$x_{\epsi_i-\epsi_j}(1)$
\end{tabular}&
\begin{tabular}[t]{ll}$s^2$,&if $\alpha=2\epsi_i$,\\
$-s$,&if $\alpha=\epsi_i+\epsi_j$,\\
$1$,&if $\alpha=2\epsi_j$.
\end{tabular}\\
\hline
5\rule{0pt}{12pt}&\small{$(A2,M5)a$}&
\begin{tabular}[t]{l}$\epsi_i+\epsi_j$, $2\epsi_k$,\\
$i<k<j$
\end{tabular}&
\begin{tabular}[t]{l}$2\epsi_i$, $\epsi_k+\epsi_j$\end{tabular}&
\begin{tabular}[t]{l}$h_{\epsi_i-\epsi_k}(s^{-1})\times$\\
$x_{\epsi_i-\epsi_k}(s)\times$\\
$x_{\epsi_k-\epsi_j}((2s^2)^{-1})\times$\\
$x_{\epsi_i-\epsi_j}(-s^{-1})$
\end{tabular}&
\begin{tabular}[t]{ll}$s^2$,&if $\alpha=2\epsi_i$,\\
$-s$,&if $\alpha=\epsi_i+\epsi_k$,\\
$1$,&if $\alpha=\epsi_i+\epsi_j$,\\
$1$,&if $\alpha=2\epsi_k$,\\
$0$,&if $\alpha=\epsi_k+\epsi_j$.
\end{tabular}\\
\hline\noalign{\eject}\hline
&&$D_{\sigma}\setminus D_{\tau}$&$D_{\tau}\setminus D_{\sigma}$&$g(s)$&$f(e_{\alpha})$, $\alpha\in\Phi^+$\\
\hline\hline
6\rule{0pt}{12pt}&\small{$(A2,M5)b$}&
\begin{tabular}[t]{l}$\epsi_i-\epsi_j$, $2\epsi_k$,\\
$i<k<j$
\end{tabular}&
\begin{tabular}[t]{l}$2\epsi_i$, $\epsi_k-\epsi_j$\end{tabular}&
\begin{tabular}[t]{l}$h_{\epsi_i+\epsi_k}(s^{-1})\times$\\
$x_{\epsi_i-\epsi_k}(s^{-1})\times$\\
$x_{\epsi_i+\epsi_j}(s^{-1})\times$\\
$h_{\epsi_k-\epsi_j}(s)$
\end{tabular}&
\begin{tabular}[t]{ll}$s^2$,&if $\alpha=2\epsi_i$,\\
$-s$,&if $\alpha=\epsi_i+\epsi_k$,\\
$1$,&if $\alpha=\epsi_i-\epsi_j$,\\
$1$,&if $\alpha=2\epsi_k$,\\
$0$,&if $\alpha=\epsi_k-\epsi_j$.
\end{tabular}\\
\hline
7\rule{0pt}{12pt}&\small{$(A2,M6)$}&
\begin{tabular}[t]{l}$\epsi_i-\epsi_k$, $2\epsi_j$,\\
$i<k<j$
\end{tabular}&
\begin{tabular}[t]{l}$2\epsi_i$\end{tabular}&
\begin{tabular}[t]{l}$h_{\epsi_i-\epsi_k}(s^{-1})\times$\\
$x_{\epsi_i+\epsi_k}(s^{-2})\times$\\
$x_{\epsi_i-\epsi_j}(1)$
\end{tabular}&
\begin{tabular}[t]{ll}$s^2$,&if $\alpha=2\epsi_i$,\\
$-s$,&if $\alpha=\epsi_i+\epsi_j$,\\
$1$,&if $\alpha=\epsi_i-\epsi_k$,\\
$1$,&if $\alpha=2\epsi_j$.
\end{tabular}\\
\hline
8\rule{0pt}{12pt}&\small{$(A3,M4)$}&
\begin{tabular}[t]{l}$\epsi_i-\epsi_k$, $2\epsi_j$,\\
$i<k<j$
\end{tabular}&
\begin{tabular}[t]{l}$\epsi_i-\epsi_j$, $2\epsi_k$\end{tabular}&
\begin{tabular}[t]{l}$h_{2\epsi_k}(s^{-2})\times$\\
$h_{\epsi_i-\epsi_j}(s^{-1})\times$\\
$x_{\epsi_k-\epsi_j}(s)$
\end{tabular}&
\begin{tabular}[t]{ll}$s^4$,&if $\alpha=2\epsi_k$,\\
$s^2$,&if $\alpha=\epsi_i-\epsi_j$,\\
$-s^2$,&if $\alpha=\epsi_k+\epsi_j$,\\
$1$,&if $\alpha=\epsi_i-\epsi_k$,\\
$1$,&if $\alpha=2\epsi_j$.
\end{tabular}\\
\hline
9\rule{0pt}{12pt}&\small{$(A3,M5)$}&
\begin{tabular}[t]{l}$\epsi_i+\epsi_j$, $2\epsi_k$,\\
$i<k<j$
\end{tabular}&
\begin{tabular}[t]{l}$\epsi_i+\epsi_k$, $2\epsi_j$\end{tabular}&
\begin{tabular}[t]{l}$h_{2\epsi_k}(s^{-1})\times$\\
$h_{2\epsi_j}(-s)\times$\\
$x_{\epsi_i-\epsi_k}(-(2s^2)^{-1})\times$\\
$x_{\epsi_i-\epsi_j}(-(2s)^{-1})\times$\\
$x_{\epsi_k-\epsi_j}(s)$
\end{tabular}&
\begin{tabular}[t]{ll}$s$,&if $\alpha=\epsi_i+\epsi_k$,\\
$1$,&if $\alpha=\epsi_i+\epsi_j$,\\
$1$,&if $\alpha=2\epsi_k$,\\
$0$,&if $\alpha=2\epsi_j$.
\end{tabular}\\
\hline
10\rule{0pt}{12pt}&\small{$(A3,M6)$}&
\begin{tabular}[t]{l}$\epsi_i-\epsi_j$, $2\epsi_k$,\\
$i<k<j$
\end{tabular}&
\begin{tabular}[t]{l}$\epsi_i+\epsi_k$\end{tabular}&
\begin{tabular}[t]{l}$x_{\epsi_i-\epsi_k}(-(2s)^{-1})\times$\\
$x_{\epsi_k+\epsi_j}(s^{-1})\times$\\
$h_{\epsi_i+\epsi_j}(s^{-1})\times$\\
$x_{\epsi_i+\epsi_j}(-1/2)$
\end{tabular}&
\begin{tabular}[t]{ll}$s$,&if $\alpha=\epsi_i+\epsi_k$,\\
$1$,&if $\alpha=\epsi_i-\epsi_j$,\\
$1$,&if $\alpha=2\epsi_k$.
\end{tabular}\\
\hline
11\rule{0pt}{12pt}&\small{$(A4,M3)$}&
\begin{tabular}[t]{l}$\epsi_i+\epsi_j$
\end{tabular}&
\begin{tabular}[t]{l}$\epsi_i+\epsi_k$,\\
$i<k<j$\end{tabular}&
\begin{tabular}[t]{l}$h_{\epsi_i-\epsi_j}(s^{-1})\times$\\
$x_{\epsi_k-\epsi_j}(-1)$
\end{tabular}&
\begin{tabular}[t]{ll}$s$,&if $\alpha=\epsi_i+\epsi_k$,\\
$1$,&if $\alpha=\epsi_i+\epsi_j$.
\end{tabular}\\
\hline
12\rule{0pt}{12pt}&\small{$(A4,M4)a$}&
\begin{tabular}[t]{l}$\epsi_k-\epsi_j$,\\
$\epsi_i+\epsi_l$,\\
$i<k<j<l$
\end{tabular}&
\begin{tabular}[t]{l}$\epsi_i+\epsi_j$,\\
$\epsi_k-\epsi_l$
\end{tabular}&
\begin{tabular}[t]{l}$x_{\epsi_j-\epsi_l}(s^{-1})\times$\\
$h_{\epsi_j-\epsi_l}(s^{-1})\times$\\
$h_{2\epsi_i}(-1)$
\end{tabular}&
\begin{tabular}[t]{ll}$s$,&if $\alpha=\epsi_k-\epsi_l$,\\
$-s$,&if $\alpha=\epsi_i+\epsi_l$,\\
$1$,&if $\alpha=\epsi_k-\epsi_j$,\\
$1$,&if $\alpha=\epsi_i+\epsi_l$.
\end{tabular}\\
\hline
13\rule{0pt}{12pt}&\small{$(A4,M4)b$}&
\begin{tabular}[t]{l}$\epsi_i-\epsi_j$,\\
$\epsi_k+\epsi_l$,\\
$i<k<j<l$
\end{tabular}&
\begin{tabular}[t]{l}$\epsi_i-\epsi_l$,\\
$\epsi_k+\epsi_j$
\end{tabular}&
\begin{tabular}[t]{l}$x_{\epsi_j-\epsi_l}(-s^{-1})\times$\\
$h_{\epsi_j-\epsi_l}(s^{-1})\times$\\
$h_{2\epsi_i}(-1)$
\end{tabular}&
\begin{tabular}[t]{ll}$s$,&if $\alpha=\epsi_k+\epsi_l$,\\
$-s$,&if $\alpha=\epsi_i-\epsi_l$,\\
$1$,&if $\alpha=\epsi_i-\epsi_j$,\\
$1$,&if $\alpha=\epsi_k+\epsi_l$.
\end{tabular}\\
\hline
14\rule{0pt}{12pt}&\small{$(A4,M4)c$}&
\begin{tabular}[t]{l}$\epsi_i+\epsi_j$,\\
$i<j$
\end{tabular}&
\begin{tabular}[t]{l}$\epsi_i+\epsi_j$
\end{tabular}&
\begin{tabular}[t]{l}$h_{2\epsi_j}(-s^{-1})\times$\\
$x_{2\epsi_j}(-s)$
\end{tabular}&
\begin{tabular}[t]{ll}$s$,&if $\alpha=\epsi_i+\epsi_j$,\\
$1$,&if $\alpha=\epsi_i-\epsi_j$.
\end{tabular}\\
\hline
15\rule{0pt}{12pt}&\small{$(A4,M5)a$}&
\begin{tabular}[t]{l}$\epsi_i+\epsi_j$,\\
$\epsi_k+\epsi_l$,\\
$i<k<j<l$
\end{tabular}&
\begin{tabular}[t]{l}$\epsi_i+\epsi_k$,\\
$\epsi_j+\epsi_l$
\end{tabular}&
\begin{tabular}[t]{l}$h_{2\epsi_i}(-1)\times$\\
$h_{\epsi_i-\epsi_l}(s^{-1})\times$\\
$h_{\epsi_k-\epsi_j}(-s^{-1})\times$\\
$x_{\epsi_k-\epsi_j}(-1)\times$\\
$h_{\epsi_i-\epsi_l}(1)$
\end{tabular}&
\begin{tabular}[t]{ll}$s^2$,&if $\alpha=\epsi_i+\epsi_k$,\\
$1$,&if $\alpha=\epsi_i+\epsi_j$,\\
$1$,&if $\alpha=\epsi_k+\epsi_l$,\\
$0$,&if $\alpha=\epsi_j+\epsi_l$.
\end{tabular}\\
\hline
16\rule{0pt}{12pt}&\small{$(A4,M5)b$}&
\begin{tabular}[t]{l}$\epsi_i+\epsi_j$,\\
$\epsi_k-\epsi_l$,\\
$i<k<j<l$
\end{tabular}&
\begin{tabular}[t]{l}$\epsi_i+\epsi_k$,\\
$\epsi_j-\epsi_l$
\end{tabular}&
\begin{tabular}[t]{l}$h_{2\epsi_i}(-1)\times$\\
$h_{\epsi_k-\epsi_j}(s^{-1})\times$\\
$x_{\epsi_i+\epsi_l}(s^{-1})\times$\\
$x_{\epsi_k-\epsi_j}(s)$
\end{tabular}&
\begin{tabular}[t]{ll}$-s$,&if $\alpha=\epsi_i+\epsi_k$,\\
$1$,&if $\alpha=\epsi_i+\epsi_j$,\\
$1$,&if $\alpha=\epsi_k-\epsi_l$,\\
$0$,&if $\alpha=\epsi_j-\epsi_l$.
\end{tabular}\\
\hline\noalign{\eject}\hline
&&$D_{\sigma}\setminus D_{\tau}$&$D_{\tau}\setminus D_{\sigma}$&$g(s)$&$f(e_{\alpha})$, $\alpha\in\Phi^+$\\
\hline\hline
17\rule{0pt}{12pt}&\small{$(A4,M6)$}&
\begin{tabular}[t]{l}$\epsi_i+\epsi_l$,\\
$\epsi_k-\epsi_j$,\\
$i<k<j<l$
\end{tabular}&
\begin{tabular}[t]{l}$\epsi_i+\epsi_k$
\end{tabular}&
\begin{tabular}[t]{l}
$h_{\epsi_k+\epsi_l}(s^{-1})\times$\\
$x_{\epsi_k-\epsi_l}(-s^{-1})\times$\\
$x_{\epsi_i+\epsi_j}(s^{-1})$
\end{tabular}&
\begin{tabular}[t]{ll}$s$,&if $\alpha=\epsi_i+\epsi_k$,\\
$1$,&if $\alpha=\epsi_i+\epsi_l$,\\
$1$,&if $\alpha=\epsi_k-\epsi_j$.
\end{tabular}\\
\hline
18\rule{0pt}{12pt}&\small{$(A5,M1)$}&
\begin{tabular}[t]{l}$\emptyset$
\end{tabular}&
\begin{tabular}[t]{l}$\epsi_i-\epsi_j$,\\
$i<j$
\end{tabular}&
\begin{tabular}[t]{l}
$h_{\epsi_i-\epsi_j}(s^{-1})$
\end{tabular}&
\begin{tabular}[t]{ll}$s$,&if $\alpha=\epsi_i-\epsi_j$.
\end{tabular}\\
\hline
19\rule{0pt}{12pt}&\small{$(A5,M2)$}&
\begin{tabular}[t]{l}$\epsi_k-\epsi_j$
\end{tabular}&
\begin{tabular}[t]{l}$\epsi_i-\epsi_j$,\\
$i<k<j$
\end{tabular}&
\begin{tabular}[t]{l}
$h_{2\epsi_i}(s^{-1})\times$\\
$x_{\epsi_i-\epsi_k}(-1)$
\end{tabular}&
\begin{tabular}[t]{ll}$s$,&if $\alpha=\epsi_i-\epsi_j$,\\
$1$,&if $\alpha=\epsi_k-\epsi_j$.
\end{tabular}\\
\hline
20\rule{0pt}{12pt}&\small{$(A5,M3)$}&
\begin{tabular}[t]{l}$\epsi_i-\epsi_k$
\end{tabular}&
\begin{tabular}[t]{l}$\epsi_i-\epsi_j$,\\
$i<k<j$
\end{tabular}&
\begin{tabular}[t]{l}
$h_{2\epsi_i}(s^{-1})\times$\\
$x_{\epsi_k-\epsi_j}(s^{-1})$
\end{tabular}&
\begin{tabular}[t]{ll}$s$,&if $\alpha=\epsi_i-\epsi_j$,\\
$1$,&if $\alpha=\epsi_i-\epsi_k$.
\end{tabular}\\
\hline
21\rule{0pt}{12pt}&\small{$(A5,M4)a$}&
\begin{tabular}[t]{l}$\epsi_i-\epsi_k$,\\
$\epsi_j-\epsi_l$,\\
$i<k<j<l$
\end{tabular}&
\begin{tabular}[t]{l}$\epsi_i-\epsi_j$,\\
$\epsi_k-\epsi_l$
\end{tabular}&
\begin{tabular}[t]{l}
$h_{2\epsi_k}(-s^{-1})\times$\\
$h_{2\epsi_i}(s^{-1})\times$\\
$x_{\epsi_k-\epsi_j}(-1)$
\end{tabular}&
\begin{tabular}[t]{ll}$s$,&if $\alpha=\epsi_i-\epsi_j$,\\
$-s$,&if $\alpha=\epsi_k-\epsi_l$,\\
$1$,&if $\alpha=\epsi_i-\epsi_k$,\\
$1$,&if $\alpha=\epsi_j-\epsi_l$.
\end{tabular}\\
\hline
22\rule{0pt}{12pt}&\small{$(A5,M4)b$}&
\begin{tabular}[t]{l}$\epsi_i-\epsi_k$,\\
$\epsi_j+\epsi_l$,\\
$i<k<j<l$
\end{tabular}&
\begin{tabular}[t]{l}$\epsi_i-\epsi_j$,\\
$\epsi_k+\epsi_l$
\end{tabular}&
\begin{tabular}[t]{l}
$h_{2\epsi_k}(-s^{-1})\times$\\
$h_{2\epsi_i}(s^{-1})\times$\\
$x_{\epsi_k-\epsi_j}(-1)$
\end{tabular}&
\begin{tabular}[t]{ll}$s$,&if $\alpha=\epsi_i-\epsi_j$,\\
$-s$,&if $\alpha=\epsi_k+\epsi_l$,\\
$1$,&if $\alpha=\epsi_i-\epsi_k$,\\
$1$,&if $\alpha=\epsi_j+\epsi_l$.
\end{tabular}\\
\hline
23\rule{0pt}{12pt}&\small{$(A5,M5)a$}&
\begin{tabular}[t]{l}$\epsi_i+\epsi_l$,\\
$\epsi_k+\epsi_j$,\\
$i<k<j<l$
\end{tabular}&
\begin{tabular}[t]{l}$\epsi_i+\epsi_j$,\\
$\epsi_k+\epsi_l$
\end{tabular}&
\begin{tabular}[t]{l}
$h_{2\epsi_k}(-s)\times$\\
$h_{2\epsi_i}(s^{-1})\times$\\
$x_{\epsi_j-\epsi_l}(-s^{-1})\times$\\
$x_{\epsi_i-\epsi_k}(s)$
\end{tabular}&
\begin{tabular}[t]{ll}$s$,&if $\alpha=\epsi_i+\epsi_j$,\\
$0$,&if $\alpha=\epsi_k+\epsi_l$,\\
$1$,&if $\alpha=\epsi_i+\epsi_l$,\\
$1$,&if $\alpha=\epsi_k+\epsi_j$.
\end{tabular}\\
\hline
24\rule{0pt}{12pt}&\small{$(A5,M5)b$}&
\begin{tabular}[t]{l}$\epsi_i-\epsi_j$,\\
$\epsi_k-\epsi_l$,\\
$i<k<j<l$
\end{tabular}&
\begin{tabular}[t]{l}$\epsi_i-\epsi_l$,\\
$\epsi_k-\epsi_j$
\end{tabular}&
\begin{tabular}[t]{l}
$h_{2\epsi_k}(-s)\times$\\
$h_{2\epsi_i}(s^{-1})\times$\\
$x_{\epsi_j-\epsi_l}(s^{-1})\times$\\
$x_{\epsi_i-\epsi_k}(s)$
\end{tabular}&
\begin{tabular}[t]{ll}$s$,&if $\alpha=\epsi_i-\epsi_l$,\\
$0$,&if $\alpha=\epsi_k-\epsi_j$,\\
$1$,&if $\alpha=\epsi_i-\epsi_j$,\\
$1$,&if $\alpha=\epsi_k-\epsi_l$.
\end{tabular}\\
\hline
25\rule{0pt}{12pt}&\small{$(A5,M5)c$}&
\begin{tabular}[t]{l}$\epsi_i-\epsi_l$,\\
$\epsi_k+\epsi_j$,\\
$i<k<j<l$
\end{tabular}&
\begin{tabular}[t]{l}$\epsi_i+\epsi_j$,\\
$\epsi_k-\epsi_l$
\end{tabular}&
\begin{tabular}[t]{l}
$h_{2\epsi_k}(-s)\times$\\
$h_{2\epsi_i}(s^{-1})\times$\\
$x_{\epsi_j+\epsi_l}(s^{-1})\times$\\
$x_{\epsi_i-\epsi_k}(s)$
\end{tabular}&
\begin{tabular}[t]{ll}$s$,&if $\alpha=\epsi_i+\epsi_j$,\\
$0$,&if $\alpha=\epsi_k-\epsi_l$,\\
$1$,&if $\alpha=\epsi_i-\epsi_l$,\\
$1$,&if $\alpha=\epsi_k+\epsi_j$.
\end{tabular}\\
\hline
26\rule{0pt}{12pt}&\small{$(A5,M5)d$}&
\begin{tabular}[t]{l}$\epsi_i-\epsi_j$,\\
$\epsi_k+\epsi_l$,\\
$i<k<j<l$
\end{tabular}&
\begin{tabular}[t]{l}$\epsi_i-\epsi_j$,\\
$\epsi_k+\epsi_l$
\end{tabular}&
\begin{tabular}[t]{l}
$h_{2\epsi_k}(-s)\times$\\
$h_{2\epsi_i}(s^{-1})\times$\\
$x_{\epsi_j+\epsi_l}(s^{-1})\times$\\
$x_{\epsi_i-\epsi_k}(s)$
\end{tabular}&
\begin{tabular}[t]{ll}$s$,&if $\alpha=\epsi_i+\epsi_l$,\\
$0$,&if $\alpha=\epsi_k-\epsi_j$,\\
$1$,&if $\alpha=\epsi_i-\epsi_j$,\\
$1$,&if $\alpha=\epsi_k+\epsi_l$.
\end{tabular}\\
\hline
27\rule{0pt}{12pt}&\small{$(A5,M6)a$}&
\begin{tabular}[t]{l}$\epsi_i-\epsi_k$,\\
$\epsi_j-\epsi_l$,\\
$i<k<j<l$
\end{tabular}&
\begin{tabular}[t]{l}$\epsi_i-\epsi_l$
\end{tabular}&
\begin{tabular}[t]{l}
$h_{2\epsi_i}(s^{-1})\times$\\
$x_{\epsi_k-\epsi_l}(s^{-1})\times$\\
$x_{\epsi_i-\epsi_j}(-1)$
\end{tabular}&
\begin{tabular}[t]{ll}$s$,&if $\alpha=\epsi_i-\epsi_l$,\\
$1$,&if $\alpha=\epsi_i-\epsi_k$,\\
$1$,&if $\alpha=\epsi_j-\epsi_l$.
\end{tabular}\\
\hline
28\rule{0pt}{12pt}&\small{$(A5,M6)b$}&
\begin{tabular}[t]{l}$\epsi_i-\epsi_k$,\\
$\epsi_j+\epsi_l$,\\
$i<k<j<l$
\end{tabular}&
\begin{tabular}[t]{l}$\epsi_i+\epsi_l$
\end{tabular}&
\begin{tabular}[t]{l}
$h_{2\epsi_i}(s^{-1})\times$\\
$x_{\epsi_k+\epsi_l}(s^{-1})\times$\\
$x_{\epsi_i-\epsi_j}(-1)$
\end{tabular}&
\begin{tabular}[t]{ll}$s$,&if $\alpha=\epsi_i+\epsi_l$,\\
$1$,&if $\alpha=\epsi_i-\epsi_k$,\\
$1$,&if $\alpha=\epsi_j+\epsi_l$.
\end{tabular}\\
\hline
29\rule{0pt}{12pt}&\small{$(A6,M1)$}&
\begin{tabular}[t]{l}$2\epsi_j$
\end{tabular}&
\begin{tabular}[t]{l}$\epsi_i+\epsi_j$,\\
$i<j$
\end{tabular}&
\begin{tabular}[t]{l}
$h_{2\epsi_i}(s^{-1})\times$\\
$x_{\epsi_i-\epsi_j}(-1/2)$
\end{tabular}&
\begin{tabular}[t]{ll}$s$,&if $\alpha=\epsi_i+\epsi_j$,\\
$1$,&if $\alpha=2\epsi_j$.
\end{tabular}\\
\hline
30\rule{0pt}{12pt}&\small{$(A6,M2)$}&
\begin{tabular}[t]{l}$2\epsi_k,2\epsi_j$
\end{tabular}&
\begin{tabular}[t]{l}$\epsi_i+\epsi_k$,\\
$i<k<j$
\end{tabular}&
\begin{tabular}[t]{l}
$h_{2\epsi_i}(s^{-1})\times$\\
$x_{\epsi_i-\epsi_k}(-1/2)\times$\\
$x_{\epsi_i-\epsi_j}(-I/2)\times$\\
$x_{\epsi_k-\epsi_j}(I)$
\end{tabular}&
\begin{tabular}[t]{ll}$s$,&if $\alpha=\epsi_i+\epsi_k$,\\
$-Is$,&if $\alpha=\epsi_i+\epsi_j$,\\
$1$,&if $\alpha=2\epsi_k$,\\
$1$,&if $\alpha=2\epsi_j$.
\end{tabular}\\
\hline\noalign{\eject}\hline
&&$D_{\sigma}\setminus D_{\tau}$&$D_{\tau}\setminus D_{\sigma}$&$g(s)$&$f(e_{\alpha})$, $\alpha\in\Phi^+$\\
\hline\hline
31\rule{0pt}{12pt}&\small{$(A6,M3)$}&
\begin{tabular}[t]{l}$\epsi_i-\epsi_k$, $2\epsi_j$,\\
$i<k<j$
\end{tabular}&
\begin{tabular}[t]{l}$\epsi_i+\epsi_j$
\end{tabular}&
\begin{tabular}[t]{l}
$h_{2\epsi_i}(s^{-1})\times$\\
$x_{\epsi_k+\epsi_j}(s^{-1})\times$\\
$x_{\epsi_i-\epsi_j}(-1/2)$
\end{tabular}&
\begin{tabular}[t]{ll}$s$,&if $\alpha=\epsi_i+\epsi_j$,\\
$1$,&if $\alpha=\epsi_i-\epsi_k$,\\
$1$,&if $\alpha=2\epsi_j$.
\end{tabular}\\
\hline
32\rule{0pt}{12pt}&\small{$(A6,M4)$}&
\begin{tabular}[t]{l}$\epsi_i-\epsi_k$,\\
$2\epsi_j$, $2\epsi_l$,\\
$i<k<j<l$
\end{tabular}&
\begin{tabular}[t]{l}$\epsi_i-\epsi_l$,\\
$\epsi_k+\epsi_j$
\end{tabular}&
\begin{tabular}[t]{l}
$h_{2\epsi_i}(s^{-1})\times$\\
$h_{2\epsi_k}(Is^{-1})\times$\\
$x_{\epsi_k-\epsi_j}(-1/2)\times$\\
$x_{\epsi_k-\epsi_l}(-I/2)\times$\\
$x_{\epsi_j-\epsi_l}(I)$
\end{tabular}&
\begin{tabular}[t]{ll}$s$,&if $\alpha=\epsi_i-\epsi_l$,\\
$-s$,&if $\alpha=\epsi_k+\epsi_l$,\\
$Is$,&if $\alpha=\epsi_i-\epsi_j$,\\
$-Is$,&if $\alpha=\epsi_k+\epsi_j$,\\
$1$,&if $\alpha=\epsi_i-\epsi_k$,\\
$1$,&if $\alpha=2\epsi_j$,\\
$1$,&if $\alpha=2\epsi_l$.
\end{tabular}\\
\hline
33\rule{0pt}{12pt}&\small{$(A6,M5)a$}&
\begin{tabular}[t]{l}$\epsi_i+\epsi_l$,\\
$2\epsi_k$, $2\epsi_j$,\\
$i<k<j<l$
\end{tabular}&
\begin{tabular}[t]{l}$\epsi_i+\epsi_k$,\\
$\epsi_j+\epsi_l$
\end{tabular}&
\begin{tabular}[t]{l}
$h_{2\epsi_i}(s^{-1})\times$\\
$h_{2\epsi_l}(Is)\times$\\
$x_{\epsi_i-\epsi_k}(-1/2)\times$\\
$x_{\epsi_k-\epsi_j}(I)\times$\\
$x_{\epsi_j-\epsi_l}(1)\times$\\
$x_{\epsi_i-\epsi_j}(-I/2)\times$\\
$x_{\epsi_k-\epsi_l}(-I)\times$\\
$x_{\epsi_j-\epsi_l}(1/2)$
\end{tabular}&
\begin{tabular}[t]{ll}$s$,&if $\alpha=\epsi_i+\epsi_k$,\\
$0$,&if $\alpha=\epsi_j+\epsi_l$,\\
$Is$,&if $\alpha=\epsi_i+\epsi_j$,\\
$1$,&if $\alpha=\epsi_i+\epsi_l$,\\
$1$,&if $\alpha=2\epsi_k$,\\
$1$,&if $\alpha=2\epsi_j$.
\end{tabular}\\
\hline
34\rule{0pt}{12pt}&\small{$(A6,M5)b$}&
\begin{tabular}[t]{l}$\epsi_i-\epsi_l$,\\
$2\epsi_k$, $2\epsi_j$,\\
$i<k<j<l$
\end{tabular}&
\begin{tabular}[t]{l}$\epsi_i+\epsi_k$,\\
$\epsi_j-\epsi_l$
\end{tabular}&
\begin{tabular}[t]{l}
$h_{2\epsi_i}(s^{-1})\times$\\
$h_{2\epsi_l}(-Is^{-1})\times$\\
$x_{\epsi_i-\epsi_k}(-1/2)\times$\\
$x_{\epsi_k-\epsi_j}(I)\times$\\
$x_{\epsi_j+\epsi_l}(-1)\times$\\
$x_{\epsi_i-\epsi_j}(-I/2)\times$\\
$x_{\epsi_k+\epsi_l}(I)$
\end{tabular}&
\begin{tabular}[t]{ll}$s$,&if $\alpha=\epsi_i+\epsi_k$,\\
$0$,&if $\alpha=\epsi_j-\epsi_l$,\\
$-Is$,&if $\alpha=\epsi_i+\epsi_j$,\\
$1$,&if $\alpha=\epsi_i-\epsi_l$,\\
$1$,&if $\alpha=2\epsi_k$,\\
$1$,&if $\alpha=2\epsi_j$.
\end{tabular}\\
\hline
35\rule{0pt}{12pt}&\small{$(A6,M6)a$}&
\begin{tabular}[t]{l}$\epsi_i-\epsi_j$,\\
$2\epsi_k$, $2\epsi_l$,\\
$i<k<j<l$
\end{tabular}&
\begin{tabular}[t]{l}$\epsi_i+\epsi_k$
\end{tabular}&
\begin{tabular}[t]{l}
$h_{2\epsi_i}(s^{-1})\times$\\
$x_{\epsi_i-\epsi_k}(-1/2)\times$\\
$x_{\epsi_k-\epsi_l}(I)\times$\\
$x_{\epsi_i-\epsi_l}(-I/2)\times$\\
$x_{\epsi_k+\epsi_j}\left(\dfrac{1-2s}{2s^2}\right)\times$\\
$x_{\epsi_i-\epsi_k}(s)\times$\\
$x_{\epsi_k+\epsi_j}\left(\dfrac{4s-1}{2s^2}\right)$\\
\end{tabular}&
\begin{tabular}[t]{ll}$s$,&if $\alpha=\epsi_i+\epsi_k$,\\
$2Is$,&if $\alpha=\epsi_k+\epsi_l$,\\
$-Is$,&if $\alpha=\epsi_i+\epsi_l$,\\
$1$,&if $\alpha=\epsi_i-\epsi_k$,\\
$1-2s$,&if $\alpha=2\epsi_j$,\\
$1+2s$,&if $\alpha=2\epsi_l$.
\end{tabular}\\
\hline
36\rule{0pt}{12pt}&\small{$(A6,M6)b$}&
\begin{tabular}[t]{l}$\epsi_i-\epsi_k$,\\
$2\epsi_j$, $2\epsi_l$,\\
$i<k<j<l$
\end{tabular}&
\begin{tabular}[t]{l}$\epsi_i+\epsi_j$
\end{tabular}&
\begin{tabular}[t]{l}
$h_{2\epsi_i}(s^{-1})\times$\\
$x_{\epsi_i-\epsi_j}(-1/2)\times$\\
$x_{\epsi_j-\epsi_l}(I)\times$\\
$x_{\epsi_i-\epsi_l}(-I/2)\times$\\
$x_{\epsi_k+\epsi_j}\left(\dfrac{1-2s}{2s^2}\right)\times$\\
$x_{\epsi_i-\epsi_j}(s)\times$\\
$x_{\epsi_k+\epsi_j}\left(\dfrac{4s-1}{2s^2}\right)$\\
\end{tabular}&
\begin{tabular}[t]{ll}$s$,&if $\alpha=\epsi_i+\epsi_j$,\\
$2Is$,&if $\alpha=\epsi_j+\epsi_l$,\\
$-Is$,&if $\alpha=\epsi_i+\epsi_l$,\\
$1$,&if $\alpha=\epsi_i-\epsi_j$,\\
$1-2s$,&if $\alpha=2\epsi_k$,\\
$1+2s$,&if $\alpha=2\epsi_l$.
\end{tabular}\\
\hline
\end{longtable}


\newpage\begin{thebibliography}{XXXX}

\bibitem[AN]{AndreNeto1} C.A.M Andr\`e., A.M. Neto. Super-characters of finite unipotent groups of types $B_n$, $C_n$ and~$D_n$. J. Algebra \textbf{305} (2006), 394--429.

\bibitem[BC]{BagnoCherniavsky1} E. Bagno, Y. Cherniavsky. Congruence $B$-orbits and the Bruhat poset of involutions of the symmetric group, see arXiv: \texttt{math.CO/0912.1819} (2009).

\bibitem[BL]{BileyLakshmibai1} S. Billey, V. Lakshmibai. Singular loci of Schubert varieties. Progr. in Math. \textbf{182}, Birkh\"auser, 2000.

\bibitem[Bo]{Bourbaki1} N. Bourbaki. Lie groups and Lie algebras. Chapters 4--6. Springer, 2002.

\bibitem[Hu]{Humphreys1} J.E. Humphreys. Linear algebraic groups. Grad. Texts in Math. \textbf{21}, Springer, 1998.

\bibitem[Ig1]{Ignatev1} M.V. Ignat'ev. Orthogonal subsets of classical root systems and coadjoint orbits of unipotent groups. Math. Notes \textbf{86} (2009), no. 1, 65--80, see also arXiv: \texttt{math.RT/0904.2841}.

\bibitem[Ig2]{Ignatev2} M.V. Ignatev. Orthogonal subsets of root systems and the orbit method. St. Petersburg Math.~J. \textbf{22} (2011), no. 5, 777--794, see also arXiv: \texttt{math.RT/1007.5220}.

\bibitem[Ig3]{Ignatev3} M.V. Ignatyev. Combinatorics of $B$-orbits and the Bruhat--Chevalley order on involutions. Transformation Groups, to appear, see also arXiv: \texttt{math.RT/1101.2189} (2011).

\bibitem[In1]{Incitti1} F. Incitti. Bruhat order on the involutions of classical Weyl groups. Ph.D. thesis. Dipartimento di Matematika ``Guido Castelnuovo'', Universit\`a di Roma ``La Sapienza'', 2003.

\bibitem[In2]{Incitti2} F. Incitti. The Bruhat order on the involutions of the symmetric groups. J. Alg. Combin. \textbf{20}~(2004), no. 3, 243--261.

\bibitem[Me1]{Melnikov1} A. Melnikov. $B$-orbit in solution to the
equation $X^2=0$ in triangular matrices. J. Algebra \textbf{223} (2000), 101--108.

\bibitem[Me2]{Melnikov2} A. Melnikov. Description of $B$-orbit
closures of order 2 in upper-triangular matrices. Transformation
Groups \textbf{11} (2006), no. 2, 217--247.

\bibitem[Me3]{Melnikov3} A. Melnikov. $B$-orbits
of nilpotent order 2 and link patterns, see arXiv:
\texttt{math.RT/0703371} (2007).

\bibitem[Pa]{Panov} A.N. Panov. Involutions in $S_n$ and associated coadjoint
orbits. J. Math. Sci. \textbf{151} (2008), no.~3, 3018--3031, see also arXiv:
\texttt{math.RT/0801.3022}.

\bibitem[RS]{RichardsonSpringer1} R.W. Richardson, T.A. Springer. The Bruhat order on symmetric varieties. Geom. Dedicata \textbf{35} (1990), no. 1--3, 389--436.


\end{thebibliography}
\end{document}